\documentclass[12pt]{article}
\usepackage{amssymb}
\usepackage{amsthm}
\usepackage{amsmath}
\usepackage{graphicx}
 \newtheorem{thm}{Theorem}[section]

 \theoremstyle{definition}
 
 \newtheorem{rem}[thm]{Remark}
 \numberwithin{equation}{section}

\begin{document}
\title{Structure of Singularities of 3D Axi-symmetric Navier-Stokes Equations}
\author{Zhen Lei\footnote{Shanghai Key Laboratory for
Contemporary Applied Mathematics; Key Laboratory of Nonlinear
Mathematical Models and Methods of Ministry of Education; School
of Mathematical Sciences, Fudan University, Shanghai 200433,
China. {\it Email:
 leizhn@yahoo.com, zlei@fudan.edu.cn}}\and
 Qi S. Zhang\footnote{Department of Mathematics, University of
 California, Riverside, CA 92521, USA. {\it Email: qizhang@math.ucr.edu}}}
\date{\today}
\maketitle

\begin{abstract}
Let $v$ be a solution of the axially symmetric Navier-Stokes
equation. We determine the structure of certain (possible) maximal
singularity of  $v$  in the following sense.   Let $(x_0, t_0)$ be
a point where the flow speed $Q_0 = |v(x_0, t_0)|$ is comparable
with the maximum flow speed at and before time $t_0$. We show
after a space-time scaling with the factor $Q_0$ and the center
$(x_0, t_0)$, the solution is arbitrarily close in $C^{2, 1,
\alpha}_{{\rm local}}$ norm to a nonzero constant vector in a
fixed parabolic cube, provided that $r_0 Q_0$ is sufficiently
large. Here $r_0$ is the distance from $x_0$ to the $z$ axis.
Similar results are also shown to be valid if $|r_0v(x_0, t_0)|$
is comparable with the maximum of $|rv(x, t)|$ at and before time
$t_0$.
\end{abstract}

\maketitle

\section{Introduction}

In this paper we study the structure, in a space time region with
maximum flow speed, of solutions to the three dimensional
incompressible Navier-Stokes equations
\begin{equation}\label{NS}
\begin{cases}
\partial_t v + v\cdot \nabla v + \nabla p = \mu \Delta v,\\[-4mm]\\
    \nabla \cdot v = 0,
\end{cases}
\quad t \geq 0,\quad x \in \mathbb{R}^3
\end{equation}
with the axially symmetric initial data
\begin{equation}\label{data}
a(x) = a^r(r, z, t)e_r +  a^\theta(r, z, t)e_\theta +  a^z(r, z,
t)e_z.
\end{equation}
In cylindrical coordinates, the solution $v=v(x, t)$ is of the
form
\begin{equation}\label{axi-solu}
v(x, t) = v^r(r, z, t)e_r +  v^\theta(r, z, t)e_\theta +  v^z(r,
z, t)e_z.
\end{equation}
Here $x = (x_1, x_2, z)$, $r = \sqrt{x_1^2 + x_2^2}$ and
\begin{equation}
e_r = \begin{pmatrix}\frac{x_1}{r}\\ \frac{x_2}{r}\\ 0
  \end{pmatrix},
\quad e_\theta = \begin{pmatrix}- \frac{x_2}{r}\\ \frac{x_1}{r}\\
0
  \end{pmatrix},
\quad e_z = \begin{pmatrix}0\\ 0\\ z
  \end{pmatrix}
\end{equation}
are the three orthogonal unit vectors along the radial, the
angular, and the axial directions respectively. Moreover, the
angular, swirl and axial components $v^r$, $v^\theta$ and $v^z$ of
the velocity field are solutions of ASNS i.e. the axially
symmetric Navier-Stokes equations
\begin{equation}\label{axi-NS}
\begin{cases}
\partial_tv^r + b\cdot\nabla v^r - \frac{(v^\theta)^2}{r}
  + \partial_rp = \big(\Delta - \frac{1}{r^2}\big)v^r,\\[-4mm]\\
\partial_tv^\theta + b\cdot\nabla v^\theta + \frac{v^rv^\theta}{r}
  = \big(\Delta - \frac{1}{r^2}\big)v^\theta,\\[-4mm]\\
\partial_tv^z + b\cdot\nabla v^z + \partial_zp = \Delta v^z,\\[-4mm]\\
b = v^re_r + v^ze_z,\quad \nabla\cdot b = \partial_rv^r +
\frac{v^r }{r} + \partial_zv^z = 0.
\end{cases}
\end{equation}
Here without loss of generality, we set the viscosity constant
$\mu = 1$.

The axially symmetric case appears much special than the full
Navier-Stokes equations, however the main regularity problem is
just as wide open. Let us briefly discuss some recent interesting
results on the axially symmetric Navier-Stokes equations. When
$v^\theta=0$, i.e. in the no swirl case, O. A. Ladyzhenskaya
\cite{L}, and M. R. Uchoviskii \& B. I. Yudovich \cite{UY} proved
that weak solutions are regular for all time. See also the work by
S. Leonardi, J. Malek, J. Necas, \& M. Pokorny \cite{LMNP}. More
recent activities, in the presence of swirl, include the results
of C.-C. Chen, R. M. Strain, T.-P.Tsai, \& H.-T. Yau in
\cite{CSTY1} \& \cite{CSTY2}, where they prove that suitable
axially symmetric solutions bounded by $Cr^{- \alpha}\sqrt{|t|}^{-
1 + \alpha}$ ($0 \leq \alpha \leq 1$) are smooth. Here $r$ is the
distance from a point to the $z$ axis, and $t$ is time. See also
the work of G. Koch, N. Nadirashvili, G. Seregin, \& V. Sverak
\cite{KNSS} and its local version in G. Seregin \& V. Sverak
\cite{SS} by different methods. Also in the presence of swirl,
there is the paper by J. Neustupa \& M. Pokorny \cite{NP}, proving
the regularity of one component (either $v^r$ or $v^{\theta}$)
implies regularity of the other components of the solution. Also
proving regularity is the work of Q. Jiu \& Z. Xin \cite{JX} under
an assumption of sufficiently small zero dimension scaled norms.
We would also like to mention the regularity results of D. Chae \&
J. Lee \cite{CL} who prove regularity results assuming finiteness
of another zero dimensional integral. On the other hand, G. Tian
\& Z. Xin \cite{TX} constructed a family of singular axis
symmetric solutions with singular initial data;  T. Hou \& C. Li
\cite{HL} found a special class of global smooth solutions. See
also a recent extension: T. Hou, Z. Lei \& C. Li \cite{HLL}.

In this paper, we take another approach to ASNS, which aims at the
understanding of the local structure of solutions when the flow
velocity is very high. This approach is more akin to the one taken
by Hamilton and Perelman in the study of Ricci flow. We are able
to do so when the flow speed $|v(x_0, t_0)|$ at a space time point
$(x_0, t_0)$ is comparable with the maximum flow speed, or
$r_0|v(x_0, t_0)|$ at a space time point $(x_0, t_0)$ is
comparable with the maximum of $r|v(x, t)|$, at and before time
$t_0$.

In order to present the result, we introduce some notations. Let
$v=v(x, t)$ be a solution to ASNS which is used here and later to
denote axially symmetric Navier-Stokes equations. Here $(x, t)$ is
a point in space time.  Given a number $a>0$, and $(x_0, t_0)$ be
a point in space time, we use the following symbol to denote the
parabolic cube
\begin{equation}\nonumber
P(x_0, t_0, a)  \equiv \{  (x, t) \ | \quad |x_0-x|< a,
 \ t_0- a^2  \le t \le t_0 \}.
\end{equation}
Unless stated otherwise,  we use $r, r_0, r_k$ to denote the
distance between points $x, x_0, x_k$ in space and the $z$ axis
respectively.

 Now we are ready to state
the main result of the paper.

\begin{thm}\label{thsingstruc}
Let $v=v(x, t)$, $(x, t) \in \mathbb{R}^3 \times [0, T_0)$, $T_0
> 0$ be a smooth  solution to the three-dimensional ASNS whose
initial condition $v_0$ satisfies
\begin{equation}\label{N0}
\|v_0\|_{L^\infty(\mathbb{R}^3)} \le N_0, \quad
\|v_0\|_{L^2(\mathbb{R}^3)} \le N_0, \quad |r v_0| \le N_0.
\end{equation}
Here $N_0$ is any given positive number. For any sufficiently
small constant $\epsilon>0$ and another constant $\sigma_0$, there
exists some $\rho_0 = \rho_0(\epsilon, N_0, \sigma_0)>0$ with the
following properties.

(a). Suppose
\begin{equation}\nonumber
r_0 |v(x_0, t_0)| \ge \rho^{-2}_0
\end{equation}
at some point
$(x_0, t_0)$ where $x_0 \in \mathbb{R}^3$ and $t_0 \in (0, T_0)$.
Suppose also $(x_0, t_0)$ is an almost maximal point in the sense:
\begin{equation}\nonumber
|v(x_0, t_0)| \ge \frac{1}{4} \sup_{x \in \mathbb{R}^3, t \le t_0}
|v(x, t)|.
\end{equation}

Then the velocity $v$ in the cube
\begin{equation}\nonumber
P(x_0, t_0, (\sigma_0 \epsilon Q)^{-1}),  \qquad Q \equiv |v(x_0,
t_0)|,
\end{equation}
is, after scaling by the factor $Q$,  $\epsilon$ close in $C^{2,1,
\alpha}_{{\rm local}}$ norm to a nonzero constant vector.

(b). The conclusion in (a) still holds if
\begin{equation}\nonumber
r_0 |v(x_0, t_0)| \ge \rho^{-2}_0
\end{equation}
at $(x_0, t_0)$ , and
\begin{equation}\nonumber
r_0 |v(x_0, t_0)| \ge \frac{1}{4} \sup_{x \in \mathbb{R}^3, t \le
t_0} r |v(x, t)|.
\end{equation}
\end{thm}

\begin{rem}
According to \cite{KNSS}, if a smooth solution blows up
in finite time, then the scaling invariant quantity $ r |v(x, t)|$
must also blow up in finite time. So, the condition in (b) can
always be satisfied.
\end{rem}

\begin{rem}
The factor $1/4$ in the statement of the theorem can be replaced
by any fixed positive number smaller than $1$.

An important open question is to generalize the current result in
$(a)$ to the case when $|v(x_0, t_0)|$ is very large but still
much smaller than maximum.

Another question is: what happens when $r_0 |v(x_0, t_0)|$ is not
large, but $|v(x_0, t_0)|$ is large at almost maximal point $(x_0,
t_0)$?
\end{rem}

\begin{rem}
The result and parameters in the theorem depend only on the
norms of the initial value in (\ref{N0}). They do not depend on
individual solutions.
\end{rem}

\section{Proof of Theorem \ref{thsingstruc}}

Let us prove part (a) first, after which the proof of (b) follows
easily.

\proof

From the condition
\begin{equation}\nonumber
\|v_0\|_{L^\infty(\mathbb{R}^3)} \le N_0, \quad
\|v_0\|_{L^2(\mathbb{R}^3)} \le N_0, \quad |r v_0 | \le N_0,
\end{equation}
by standard theory (see Proposition 4.1 in \cite{KNSS}
e.g.), there exists a time $h_0$ such that
\begin{equation}\label{vt2C0}
\|v(\cdot, t)\|_{L^\infty(\mathbb{R}^3)} \le 2 N_0, \qquad t \le
h_0.
\end{equation}

The proof is divided into several steps, using the method of
contradiction.
\medskip

{\it step 1.}  setting up a limit solution.

Suppose part (a) of the theorem is false. Then for some $\epsilon
> 0$ and $\sigma_0 > 0$, there exists a sequence of solutions $v_k$ with
initial condition satisfying \eqref{N0}, defined on the time
 interval $[0, T_k)$ for some $T_k > h_0$,
 which satisfies the following conditions.

(i) there exist sequences of positive numbers $\rho_k \to 0$ ,
points $x_k \in \mathbb{R}^3$, and times $t_k \in [0, T_k)$ such
that
\begin{equation}\nonumber
r_k |v_k(x_k, t_k)| \ge \rho^{-2}_k;
\end{equation}

 (ii).  for each $k$,  the solution $v_k$ in the parabolic region
\begin{equation}\nonumber
P(x_k, t_k, [c Q_k]^{-1}) \equiv \{  (x, t) \in [0, T_k) \ |   \
|x_k-x|< (c Q_k)^{-1}, \  t_k- (c Q_k)^{-2} \le t \le t_k \}
\end{equation}
is not, after scaling by the factor $Q_k$, $\epsilon$  close,
in $C^{2, 1, \alpha}$ norm, to a nonzero constant vector.  Here $c
= \sigma_0 \epsilon$ and also
\begin{equation}\nonumber
Q_k=|v_k(x_k, t_k)| \ge \frac{1}{4} \sup_{t \in [0, t_k], \ x \in
R^3}  |v_k(x, t)|.
\end{equation}

Write $\alpha_k = r_k Q_k = r_k |v_k(x_k, t_k)|$.  We consider
$v_k$ in the space time cube
\begin{equation}\nonumber
 P(x_k, t_k, r_k/\sqrt{\alpha_k}) \equiv  B(x_k,  r_k/\sqrt{\alpha_k}) \times
 [t_k- (r_k/\sqrt{\alpha_k})^2, t_k].
\end{equation}
Note that
 \begin{equation} \label{betakQk}
 \beta_k \equiv \frac{r_k}{\sqrt{\alpha_k}} = \frac{r_k}{\sqrt{r_k Q_k}} = o(r_k),
 \end{equation}
\begin{equation}\nonumber
Q_k \beta_k = \sqrt{r_k Q_k} \to \infty,  \quad k \to \infty.
\end{equation}
Define the scaled function
\begin{equation}\label{vktild}
\tilde v_k = Q_k^{-1} v_k(Q_k^{-1} \tilde x + x_k,
Q_k^{-2} \tilde t + t_k)
\end{equation}
Then $\tilde v_k$ is solution of the Navier-Stokes equation in the
slab $\mathbb{R}^3 \times [ -( Q_k \beta_k)^2, 0]$. Moreover,  by
the assumption on $Q_k$, we know that $|\tilde v_k| \le 4$
whenever defined. Since $\tilde v_k$ are bounded mild solutions,
we know from Proposition 4.1 in \cite{KNSS} e.g. that the $C^{2,1,
\alpha}$ norm of $v_k$ are uniformly bounded in $\mathbb{R}^3
\times [ -( Q_k \beta_k)^2 +1, 0]$.  In addition, the pressure
$P_k$, satisfying $\Delta P_k = div (v_k \nabla_k)$, also has
uniformly bounded $C^{2,1, \alpha}_{{\rm local}}$ norm, by virtue
of standard Schauder theory.   Actually all $C^{p, p/2}$ norms are
bounded for $p \ge 1$. But we do not need this fact here.

Let us restrict the solution $\tilde v_k$ to the cube
\begin{equation}\nonumber
P(0, 0,  Q_k \beta_k) = \{ (\tilde x, \tilde t) \ | \  |\tilde x|
\le Q_k \beta_k, -( Q_k \beta_k)^2 \le \tilde t \le 0 \}.
\end{equation}
By the uniform bounds on $C^{2,1, \alpha}_{{\rm local}}$
norm and the fact that $Q_k \beta_k \to \infty$, we know there
exists a subsequence, still called $\{ \tilde v_k \}$ that
converges to an ancient solution of Navier-Stokes equation in
$C^{2,1, \alpha}_{{\rm local}}$ sense. Let us call this ancient
solution $\tilde v$. In the next step, we will show that it is 2
spatial dimensional solution, with one of the dimension being the
$z$ axis.
\medskip

{\it step 2.}  proving $\tilde v$ is a 2D solution.
\medskip

Denote by $v^\theta_k$ the angular component of $v_k$. For the
given initial value, it is known that
\begin{equation}\nonumber
|v^\theta_k(x, t)| \le \frac{ N_0}{r}.
\end{equation}
For $x \in
B(x_k, \beta_k)$, we have, by (\ref{betakQk}),
\begin{equation}\nonumber
|v^\theta_k(x, t)| \le \frac{ 2 N_0}{r_k}
\end{equation}
when $k$
is sufficiently large. Therefore
\begin{equation}\label{Qvtheta}
Q^{-1}_k |v^\theta_k(x, t)| \le \frac{ 2 N_0}{Q_k
r_k} \to 0, \qquad k \to \infty.
\end{equation}

In the standard basis for $\mathbb{R}^3$, let $x_k=(x_{k, 1},
x_{k, 2}, x_{k, 3})$ with the third component being the one for
the $z$ axis, and let $\xi_k = (0, 0, x_{k, 3})$. Since the
vectors $(x_k-\xi_k)/|x_k-\xi_k|$  are unit ones, there exists a
subsequence, still labeled by $k$, which converges to a unit
vector $\zeta=(\zeta_1, \zeta_2, 0)$.  We use
\begin{equation}\nonumber
\zeta,  \zeta'=(-\zeta_2, \zeta_1, 0), (0, 0, 1)
\end{equation}
as
the basis of a new coordinate.  Since this basis is obtained by a
rotation around $z$ axis, we know $v_k$ is invariant.  From now
on, when we mention the coordinates of a point, we mean to use the
new basis with the same origin. We still use $(\theta, r, z)$ to
denote the variables for the cylindrical system corresponding to
this new basis.

For $x \in B(x_k, \beta_k)$, we recall that $\theta$ is the angle
between $x$ and $\zeta$. Then
\begin{equation}\label{thetakto0}
\cos \theta =
 \frac{(x-\xi_k) \cdot
\zeta}{|x - (0, 0, x_3)^T|} = \frac{(x_k-\xi_k) \cdot \zeta}{|x_k
-\xi_k|} + \frac{O(\beta_k)}{r_k} \to 1, \qquad k \to \infty.
\end{equation}

For $v_k=v_k(x, t)$  in $B(x_k, \beta_k) \times [t_k-\beta^2_k,
t_k]$, we have defined
\begin{equation}\nonumber
\tilde v_k= \tilde v_k(\tilde x, \tilde t) = Q_k^{-1} v_k(Q_k^{-1}
\tilde x + x_k,  Q_k^{-2} \tilde t + t_k)
\end{equation}
where
$x=Q_k^{-1} \tilde x + x_k$ and $t=Q_k^{-2} \tilde t + t_k$. Then
for $x=(x^{(1)}, x^{(2)}, x^{(3)})$ and $\tilde x=(\tilde x^{(1)},
\tilde x^{(2)}, \tilde x^{(3)})$, we have
\begin{equation}\label{vkdaoshu}
\begin{cases}
\partial_r v_k(x, t) = \partial_{x^{(1)}} v_k(x, t)
  \cos \theta + \partial_{x^{(2)}} v_k(x, t)  \sin \theta\\
\quad\quad = Q^2_k \partial_{\tilde x^{(1)}} \tilde v_k(\tilde x,
  \tilde t) \cos \theta + Q^2_k \partial_{\tilde x^{(2)}} \tilde v_k (\tilde x, \tilde t) \sin
  \theta\\[-4mm]\\
\partial^2_r v_k(x, t) = Q^3_k \partial^2_{\tilde x^{(1)}} \tilde v_k(\tilde x, \tilde t)  \cos^2 \theta +
  2 Q^3_k \partial^2_{\tilde x^{(1)} \tilde x^{(2)}} \tilde v_k (\tilde x, \tilde t) \sin \theta
  \cos \theta\\
\quad\quad +\ Q^3_k \partial^2_{\tilde x^{(2)}} \tilde v_k(\tilde
x, \tilde t)  \sin^2
\theta;\\[-4mm]\\
\partial^2_z v_k(x, t) = Q^3_k \partial^2_{\tilde x^{(3)}} \tilde v_k(\tilde x, \tilde
t);\\[-4mm]\\
\partial^2_t v_k(x, t) = Q^3_k \partial^2_{\tilde t} \tilde v_k(\tilde x, \tilde t).
\end{cases}
\end{equation}

For the pressure $p_k=p_k(x, t)$, recall that
\begin{equation}\nonumber
\tilde p_k =
 \tilde p_k(\tilde x, \tilde t) = Q_k^{-2} p_k(Q_k^{-1} \tilde x + x_k,  Q_k^{-2} \tilde t + t_k).
\end{equation}
Therefore
\begin{equation}\label{pdaoshu}
\partial_r p_k(x, t) =Q^3_k \partial_{\tilde x^{(1)}} \tilde p_k(\tilde x, \tilde t)  \cos \theta +
Q^3_k \partial_{\tilde x^{(2)}} \tilde p_k (\tilde x, \tilde t)
\sin \theta
\end{equation}

Writing $v_k= v^r_k e_r +  v^\theta_k e_\theta + v^z_k e_z$, then
\begin{eqnarray}\label{driftdaoshu}
&&v^r_k \partial_r v^r_k + v^z_k \partial_z v^r_k = Q^3_k [
  v^r_k(\tilde x, \tilde t) \partial_{\tilde x^{(1)}} \tilde
  v^r_k(\tilde x, \tilde t)  \cos \theta\\\nonumber
&&\quad +\ v^r_k(\tilde x, \tilde t)
  \partial_{\tilde x^{(2)}} \tilde v^r_k (\tilde x, \tilde t) \sin
  \theta + \tilde v^z_k \partial_{\tilde x^{(3)}} \tilde
  v^r_k(\tilde x, \tilde t) ].
\end{eqnarray}
We substitute the above identities into the equation for $v^r_k$:
\begin{equation}\nonumber
[\partial^2_r + (1/r) \partial_r-\frac{1}{r^2}]
v^r_k-(b\cdot\nabla)v^r_k+
\frac{(v^{\theta}_k)^2}{r}-\frac{\partial p_k}{\partial
r}-\frac{\partial v^r_k}{\partial t}=0,
\end{equation}
we arrive at
\begin{eqnarray}\nonumber
&&(\partial^2_{\tilde x^{(1)}} + \partial^2_{\tilde x^{(3)}} )
  \tilde v^r_k -( \tilde v^r_k \partial_{\tilde x^{(1)}} + \tilde
  v^z_k \partial_{\tilde x^{(3)}} ) \tilde v^r_k -\partial_{\tilde
  x^{(1)} }\tilde p_k - \partial_{\tilde t} v^r_k + \frac{1}{Q_k r}( \partial_{\tilde x^{(1)}} \tilde
  v_k(\tilde x, \tilde t)  \cos \theta \\\nonumber
&&\quad +\ \partial_{\tilde x^{(2)}}
  \tilde v_k (\tilde x, \tilde t) \sin \theta) - \frac{1}{(Q_k r)^2}
  \tilde v^r_k+ \frac{(r v^\theta_k)^2}{(Q_k r)^3}+ O(\theta)=0.
\end{eqnarray}
Here the term $O(\theta)$ stands for all those terms which
vanish when $\theta \to 0$ as $k \to \infty$.  In particular all
terms involving the derivative with respect to $\tilde x^{(2)}$
are included in $O(\theta)$.

Recall that $Q_k r $ is comparable to $Q_k r_k$ which goes to
$\infty$. Letting $k \to \infty$ and noting that $v_k$ and
derivatives are uniformly bounded, we know that $\tilde v^{1}$,
the limit of $\tilde v^r_k$ satisfies
\begin{equation}\nonumber
(\partial^2_{\tilde x^{(1)}} +\partial^2_{\tilde x^{(3)}} ) \tilde
v^{(1)} -( \tilde v^{(1)} \partial_{\tilde x^{(1)}} + \tilde
v^{(3)} \partial_{\tilde x^{(3)}} ) \tilde v^{(1)}
-\partial_{\tilde x^{(1)} }\tilde p -
\partial_{\tilde t} v^{(1)} =0.
\end{equation}
Here $\tilde v^{(3)}$ is the limit of $v^z_k$, for which we
have, in a similar manner
\[
(\partial^2_{\tilde x^{(1)}} +\partial^2_{\tilde x^{(3)}} ) \tilde
v^{(3)} -( \tilde v^{(1)} \partial_{\tilde x^{(1)}} + \tilde
v^{(3)} \partial_{\tilde x^{(3)}} ) \tilde v^{(3)}
-\partial_{\tilde x^{(3)} }\tilde p -
\partial_{\tilde t} v^{(3)} =0.
\]

Note that $\tilde v_k$ and its derivatives are uniformly bounded
in the region of concern. When $k \to \infty$, $\theta \to 0$ in
the region of concern.  Hence $\tilde v^\theta_k$ and derivatives
all vanish when $k \to \infty$.

Finally we need to show that $\tilde v^{(1)}$ and $\tilde v^{(3)}$
are independent of the variable $\tilde x^{(2)}$.  To prove it,
let us recall that $\partial_\theta v^r_k=\partial_\theta
v^z_k=0$. Hence
\begin{equation}\nonumber
-\partial_{x^{(1)}} v_k^r \sin \theta + \partial_{x^{(2)}} v_k^r
\cos \theta = -\partial_{x^{(1)}} v_k^z \sin \theta +
\partial_{x^{(2)}} v_k^z \cos \theta = 0.
\end{equation}
This implies
\begin{equation}\nonumber
\partial_{\tilde x^{(2)}}  \tilde v_k = \partial_{\tilde x^{(1)} } \tilde v_k \tan
\theta. \end{equation}
Taking $k \to \infty (\theta \to 0)$ we see
the desired result.

 \medskip

{\it Step 3.}

Here we just use the fact that 2 dimensional ancient (mild)
solutions are constants (\cite{KNSS}) and the regularity result
Proposition 4.1 in the same paper to conclude that $\tilde v_k$,
with $k$ large, is $\epsilon$ close to a nonzero constant vector
in $C^{2, 1, \alpha}_{{\rm local}}$ sense. This contradiction with
the condition (ii) at the beginning of the section proves part (a)
of the theorem.

\medskip

Now we prove part (b).

Suppose part (b) of the theorem is false. Then for some
$\epsilon>0$, there exists a sequence of solutions $v_k$ with
normalized initial condition as above, defined on the time
 interval $[0, T_k)$ for some $T_k \in [h_0, T_0]$,
  which satisfies the following conditions.

(i) there exist sequences of positive numbers $\rho_k \to 0$ ,
points $x_k \in \mathbb{R}^3$, and times $t_k \in [0, T_k)$ such
that
\begin{equation}\nonumber
r_k |v_k(x_k, t_k)| \ge \rho^{-2}_k;
\end{equation}

 (ii).  for each $k$,  the solution $v_k$ in the parabolic region
\begin{equation}\nonumber
P(x_k, t_k, [c Q_k]^{-1}) \equiv \{  (x, t) \in [0, T_k) \ |   \
d(x_k,  x, t_k)< (c Q_k)^{-1}, \  t_k- (c Q_k)^{-2} \le t \le t_k
\}
\end{equation}
is not, after scaling by the factor $Q_k$, $\epsilon$ close, in
$C^{2, 1, \alpha}_{{\rm local}}$ norm, to a nonzero constant
vector. Here $c = \sigma_0 \epsilon$ and also
\begin{equation}\nonumber
 r_k |v_k(x_k, t_k)| \ge \frac{1}{4} \sup_{t \in [0, t_k], \ x \in R^3}
 r |v_k(x, t)|.
\end{equation}

Define as before $Q_k =|v(x_k, t_k) |$.  Suppose  $k$ is large.
Then for $x \in B(x_k, \beta_k)$ with $\beta_k = r_k/\sqrt{ r_k
Q_k} = o(r_k)$, there holds, for $t \le t_k$,
\begin{equation}\nonumber
r |v(x, t)| \le r_k |v(x_k, t_k) | = r_k Q_k
\end{equation}
and
\begin{equation}\nonumber
r_k/2 \le r \le 2 r_k
\end{equation}
when $k$ large. This shows, in the ball $B(x_k, \beta_k)$ and for $t \le t_k$,
\begin{equation}\nonumber
|v(x, t)| \le 2 Q_k.
\end{equation}
Now we can scale by $Q^{-1}_k$ in the above ball again as in the
proof of part (a). By Theorem 2.8 of \cite{SS}, the limit of
scaled solutions is again a bounded, mild, ancient solution.
Similar arguments as in part (a) lead to a contradiction, proving
part (b).
 \qed

{\bf Acknowledgment} Qi S. Zhang would like to thank the
mathematics department of Fudan university for its hospitality
during his visits. The authors would like to thank Professor
Thomas Y. Hou for his interest. The work was in part supported by
NSFC (grants No. 10801029 and 10911120384), FANEDD, Shanghai
Rising Star Program (10QA1400300), SGST 09DZ2272900 and SRF for
ROCS, SEM.

\end{document}